\documentclass[12pt]{article}
\usepackage{graphicx}
\usepackage{amsmath}

\input{psfig.sty}
\newtheorem{theorem}{Theorem}[section]
\newtheorem{lemma}[theorem]{Lemma}

\begin{document}

\begin{center}
\begin{Large}
{\bf Eigenvalues of Sturm Liouville problems\\
with discontinuity conditions inside a finite interval}\\[0pt]
\end{Large}
\vspace*{1cm} by\\[0pt]
\vspace*{1cm} \textbf{B. Chanane}\\\vspace*{0.5cm}
Mathematical Science Department,\\[0pt]
K.F.U.P.M., Dhahran 31261, Saudi Arabia\\[0pt]
E-Mail: chanane@kfupm.edu.sa\\[0pt]
\end{center}

\pagenumbering{arabic} \setcounter{page}{1}

\textbf{Abstract} --- In this work, we use the \textit{regularized
sampling method} to compute the eigenvalues of Sturm Liouville
problems with discontinuity conditions inside a finite interval.
We work out an example by computing a few eigenvalues and their
corresponding eigenfunctions.


\section{Introduction}

\setcounter{equation}{0}

In \cite{C05a}, this author introduced the \textit{regularized
sampling method}. A method which is based on Shannon's sampling
theory but applied to regularized functions, hence avoiding any
(multiple) integration(s) and keeping the number of terms in the
cardinal series manageable. It has been demonstrated that the
method is capable of delivering higher order estimates of the
eigenvalues at a very low cost. The purpose in this paper is to
extend the domain of application of this method to the problem at
hand.

\section{Main results}

Consider the following Sturm-Liouville problem with discontinuity inside a
finite interval,

\begin{equation}
\left\{
\begin{array}{c}
-y^{\prime \prime }+q(x)y=\mu ^{2}y\qquad ,\qquad x\in (0,\pi ) \\
y^{\prime }(0)=0=y(\pi ) \\
y(d+0)=ay(d-0) \\
y^{\prime }(d+0)=a^{-1}y^{\prime }(d-0)%
\end{array}%
\right.  \label{eq1}
\end{equation}%
where $a>0$, $a\neq 1$, $0<d<\pi $ and $q\in L^{2}(0,\pi )$.

It has been shown in \cite{A2006} that the problem has a point
spectrum and each eigenvalue has multiplicity one and accumulate
only at $+\infty$. The purpose in this paper is to compute the
eigenvalues of (\ref{eq1}) with the minimum of
effort and a greater precision using the newly introduced \textit{%
regularized sampling method} \cite{C05a}, an improvement on the method based
on sampling theory introduced in \cite{BC96}.

Let $y_{L}(x,\mu )$~and $y_{R}(x,\mu )$, be the solutions of the base
problems
\begin{equation}
\left\{
\begin{array}{c}
-y_{L}^{\prime \prime }+q(x)y_{L}=\mu ^{2}y_{L}\qquad ,\qquad x\in \lbrack
0,d] \\
~y_{L}(0)=1,\text{ }y_{L}^{\prime }(0)=0%
\end{array}%
\right.  \label{eq2}
\end{equation}%
and
\begin{equation}
\left\{
\begin{array}{c}
-y_{R}^{\prime \prime }+q(x)y_{R}=\mu ^{2}y_{R}\qquad ,\qquad x\in \lbrack
d,\pi ] \\
y_{R}(\pi )=0~,~y_{R}^{\prime }(\pi )=1%
\end{array}%
\right.  \label{eq3}
\end{equation}%
respectively.

\bigskip Using $y^{\prime }(\pi ,\mu )=\alpha $, an unknown constant to be
determined along the eigenvalue parameter $\mu ^{2}$ the discontinuity
conditions give%
\begin{equation}
\begin{array}{c}
\alpha y_{R}(d,\mu )=ay_{L}(d,\mu ) \\
\alpha y_{R}^{\prime }(d,\mu )=a^{-1}y_{L}^{\prime }(d,\mu )%
\end{array}
\label{eq4}
\end{equation}%
Note that $\alpha \neq 0$ since otherwise, $y_{L}(d,\mu )=y_{L}^{\prime
}(d,\mu )=0$ leading to $y_{L}(x,\mu )\equiv 0$ which contradicts $%
y_{L}(0,\mu )=1$. \ A necessary and sufficient condition for non trivial
solutions is $\Delta (\mu )=0$ where the characteristic function $\Delta $
is defined by,%
\begin{equation}
\Delta (\mu )=\left\vert
\begin{array}{cc}
ay_{L}(d,\mu ) & y_{R}(d,\mu ) \\
a^{-1}y_{L}^{\prime }(d,\mu ) & y_{R}^{\prime }(d,\mu )%
\end{array}%
\right\vert  \label{eq5}
\end{equation}%
Thus the eigenvalues of the problem at hand are seen as the square
of the zeroes of $\Delta $. Let $\mu _{k\text{ }}$ be such a zero,
$\alpha $ will take the value
\begin{equation}
\alpha _{k}=\frac{ay_{L}(d,\mu _{k\text{ }})}{y_{R}(d,\mu _{k\text{ }})}
\label{eq6}
\end{equation}%
The eigenfunction associated to the simple eigenvalue $\mu _{k\text{ }}$is
defined by
\begin{equation}
y_{[k]}(x,\mu _{k\text{ }})=\left\{
\begin{array}{c}
y_{L}(x,\mu _{k\text{ }})\text{ if }0\leq x\leq d \\
\alpha _{k}y_{R}(x,\mu _{k\text{ }})\text{ if }d\leq x\leq \pi%
\end{array}%
\right.  \label{eq7}
\end{equation}%
In the following we shall present some important properties of the functions
$y_{L}(x,\mu )$ and $y_{R}(x,\mu )$ .

We shall need the following well known lemmata,

\begin{lemma}
$\sin z/z$~and $\cos z$~are entire as functions of $z$~and satisfy the
estimates
\begin{equation}
\left\vert \sin z/z\right\vert \leq \gamma _{0}e^{\left\vert
Im~z\right\vert }/(1+\left\vert z\right\vert
)~\mbox{and}~\left\vert \cos z\right\vert \leq e^{\left\vert
Im~z\right\vert }  \notag
\end{equation}%
where $\gamma _{0}=1.72$.\bigskip
\end{lemma}

\begin{lemma}
If $\beta $ is a positive constant and $\varphi $ is a positive
function
satisfying,%
\begin{equation*}
z(x)\leq \beta +\int_{x}^{\pi }\varphi (t)z(t)dt
\end{equation*}%
then%
\begin{equation*}
z(x)\leq \beta \exp \left\{ \int_{x}^{\pi }\varphi (t)dt\right\}
\end{equation*}
\end{lemma}

\begin{theorem}
$y_{L}(x,\mu )$, $y_{L}^{\prime }(x,\mu )$, $y_{L}(x,\mu )-\cos \mu x$, $%
y_{L}^{\prime }(x,\mu )+\mu \sin \mu x$ $\in PW_{x}$ as functions of $\mu $%
~for each fixed $x\in (0,d]$~and satisfy the growth conditions,%
\begin{eqnarray*}
\left\vert y_{L}(x,\mu )\right\vert  &\leq &\gamma
_{1}e^{x|\bf{Im}\mu |}
\\
\left\vert y_{L}(x,\mu )-\cos \mu x\right\vert  &\leq &\frac{\gamma _{2}}{%
1+|\mu |\pi }e^{x|\bf{Im}\mu |} \\
\left\vert y_{L}^{\prime }(x,\mu )+\mu \sin \mu x-\int_{0}^{x}\cos \mu
(x-t)q(t)\cos (\mu t)dt\right\vert  &\leq &\frac{\gamma _{3}}{1+|\mu |\pi }%
e^{x|\bf{Im}\mu |} \\
\left\vert y_{L}^{\prime }(x,\mu )+\mu \sin \mu x\right\vert
&\leq &\gamma _{4}e^{x|\bf{Im}\mu |}
\end{eqnarray*}%
\newline
$y_{R}(x,\mu )$, $y_{R}^{\prime }(x,\mu )+\frac{\sin \mu (\pi -x)}{\mu }$, $%
y_{R}^{\prime }(x,\mu )-\cos \mu (\pi -x)$ $\in PW_{\pi -x}$ as functions of
$\mu $~for each fixed $x\in \lbrack d,\pi )$~and satisfy the growth
conditions,%
\begin{equation*}
\left\vert y_{R}(x,\mu )\right\vert \leq \frac{\gamma _{5}}{1+|\mu |\pi }%
e^{(\pi -x)|\bf{Im}\mu |}
\end{equation*}%
\begin{equation*}
\left\vert y_{R}(x,\mu )+\frac{\sin \mu (\pi -x)}{\mu }\right\vert
\leq \frac{\gamma _{6}}{\left( 1+|\mu |\pi \right) ^{2}}e^{(\pi
-x)|\bf{Im}\mu |}
\end{equation*}%
\begin{equation*}
\left\vert y_{R}^{\prime }(x,\mu )-\cos \mu (\pi -x)\right\vert \leq \frac{%
\gamma _{7}}{1+|\mu |\pi }e^{(\pi -x)|\bf{Im}\mu |}
\end{equation*}%
where $\gamma _{1},\cdot \cdot \cdot ,\gamma _{7}$ are some positive
constants.
\end{theorem}

\bigskip
\textbf{Proof:} For $y_{L}$ we have,
\begin{equation}
y_{L}(x,\mu )=\cos \mu x+\int_{0}^{x}\frac{\sin \mu (x-t)}{\mu }%
q(t)y_{L}(t,\mu )dt  \label{eq8}
\end{equation}%
and
\begin{equation}
y_{L}^{\prime }(x,\mu )=-\mu \sin \mu x+\int_{0}^{x}\cos \mu
(x-t)q(t)y_{L}(t,\mu )dt\text{.}  \label{eq9}
\end{equation}%
Standard arguments show that $y_{L}$, $y_{L}^{\prime }$ are entire
functions of $\mu $ for each $x\in (0,d]$ and belong to $PW_{x}$ for each $%
x\in (0,d]$ as function of $\mu $.%
\begin{equation*}
\left\vert y_{L}(x,\mu )\right\vert \leq e^{x|\bf{Im}\mu
|}+\int_{0}^{x}\gamma _{0}\pi e^{(x-t)|\bf{Im}\mu |}\left\vert
q(t)\right\vert \cdot \left\vert y_{L}(t,\mu )\right\vert dt
\end{equation*}%
\begin{equation*}
e^{-x|\bf{Im}\mu |}\left\vert y_{L}(x,\mu )\right\vert \leq
1+\int_{0}^{x}\gamma _{0}\pi \left\vert q(t)\right\vert \cdot e^{-t|\bf{Im}%
\mu |}\left\vert y_{L}(t,\mu )\right\vert dt
\end{equation*}%
using Gronwall's lemma gives,%
\begin{eqnarray*}
e^{-x|\bf{Im}\mu |}\left\vert y_{L}(x,\mu )\right\vert  &\leq
&\exp
\left\{ \gamma _{0}\pi \int_{0}^{x}\left\vert q(t)\right\vert dt\right\}  \\
&\leq &\exp \left\{ \gamma _{0}\pi \int_{0}^{\pi }\left\vert q(t)\right\vert
dt\right\}
\end{eqnarray*}%
so that,%
\begin{equation}
\left\vert y_{L}(x,\mu )\right\vert \leq \exp \left\{ \gamma
_{0}\pi \int_{0}^{\pi }\left\vert q(t)\right\vert dt\right\}
e^{x|\bf{Im}\mu |}=\gamma _{1}e^{x|\bf{Im}\mu |}
\end{equation}%
where $\gamma _{1}=\exp \left\{ \gamma _{0}\pi \int_{0}^{\pi }\left\vert
q(t)\right\vert dt\right\} $.%
\begin{equation}
y_{L}(x,\mu )-\cos \mu x=\int_{0}^{x}\frac{\sin \mu (x-t)}{\mu }%
q(t)y_{L}(t,\mu )dt  \label{eq10}
\end{equation}%
\begin{eqnarray}
\left\vert y_{L}(x,\mu )-\cos \mu x\right\vert  &\leq
&\int_{0}^{x}\left\vert \frac{\sin \mu (x-t)}{\mu }\right\vert \cdot
\left\vert q(t)\right\vert \cdot \left\vert y_{L}(t,\mu )\right\vert dt
\notag \\
&\leq &\int_{0}^{x}\left\vert \frac{\sin \mu (x-t)}{\mu
(x-t)}\right\vert (x-t)\cdot \left\vert q(t)\right\vert \cdot
\gamma _{1}e^{t|\bf{Im}\mu |}dt
\\
&\leq &\int_{0}^{x}\gamma _{0}\frac{(x-t)}{1+|\mu |(x-t)}\cdot \left\vert
q(t)\right\vert \cdot \gamma _{1}e^{t|\bf{Im}\mu |}dt \\
&\leq &e^{x|\bf{Im}\mu |}\gamma _{0}\gamma _{1}\frac{\pi }{1+|\mu |\pi }%
\int_{0}^{\pi }\left\vert q(t)\right\vert dt=\frac{\gamma _{2}}{1+|\mu |\pi }%
e^{x|\bf{Im}\mu |}
\end{eqnarray}

where $\gamma _{2}=\pi \gamma _{0}\gamma _{1}\int_{0}^{\pi }\left\vert
q(t)\right\vert dt$.

Now,%
\begin{equation}
y_{L}^{\prime }(x,\mu )+\mu \sin \mu x-\int_{0}^{x}\cos \mu
(x-t)q(t)\cos (\mu t)dt=\int_{0}^{x}\cos \mu (x-t)q(t)\left\{
y_{L}(t,\mu )-\cos \mu t\right\} dt
\end{equation}%
so that,%
\begin{eqnarray}
&&\left\vert y_{L}^{\prime }(x,\mu )+\mu \sin \mu x-\int_{0}^{x}\cos \mu
(x-t)q(t)\cos (\mu t)dt\right\vert   \notag \\
&\leq &\int_{0}^{x}e^{(x-t)|\bf{Im}\mu |}\left\vert
q(t)\right\vert
\left\vert y_{L}(t,\mu )-\cos \mu t\right\vert dt  \notag \\
&\leq &e^{x|\bf{Im}\mu |}\frac{\gamma _{2}}{1+|\mu |\pi
}\int_{0}^{\pi }|q(t)|dt=\frac{\gamma _{3}}{1+|\mu |\pi
}e^{x|\bf{Im}\mu |}
\end{eqnarray}%
where $\gamma _{3}=\gamma _{2}\int_{0}^{\pi }\left\vert q(t)\right\vert dt$.%
\begin{eqnarray}
\left\vert y_{L}^{\prime }(x,\mu )+\mu \sin \mu x\right\vert  &=&\left\vert
\int_{0}^{x}\cos \mu (x-t)q(t)y_{L}(t,\mu )dt\right\vert  \\
&\leq &\gamma _{1}e^{x|\bf{Im}\mu |}\int_{0}^{\pi }\left\vert
q(t)\right\vert dt=\gamma _{4}e^{x|\bf{Im}\mu |}
\end{eqnarray}%
where $\gamma _{4}=\gamma _{1}\int_{0}^{\pi }\left\vert q(t)\right\vert dt$.

As for $y_{R}$ we have,%
\begin{equation}
y_{R}(x,\mu )=-\frac{\sin \mu (\pi -x)}{\mu }-\int_{x}^{\pi }\frac{\sin \mu
(x-t)}{\mu }q(t)y_{R}(t,\mu )dt  \label{eq11}
\end{equation}%
and standard arguments show that $y_{R}$ and $y_{R}^{\prime }$ are
entire functions of $\mu $ for each $x\in \lbrack d,\pi )$ and belong to $%
PW_{\pi -x}$ for each $x\in \lbrack d,\pi )$ as function of $\mu $. Also,
\begin{equation}
y_{R}(x,\mu )+\frac{\sin \mu (\pi -x)}{\mu }=-\int_{x}^{\pi }\frac{\sin \mu
(x-t)}{\mu }q(t)y_{R}(t,\mu )dt  \label{eq12}
\end{equation}%
and%
\begin{equation}
y_{R}^{\prime }(x,\mu )-\cos \mu (\pi -x)=\int_{x}^{\pi }\cos \mu
(x-t)q(t)y_{R}(t,\mu )dt  \label{eq13}
\end{equation}%
are entire functions of $\mu $ for each $x\in \lbrack d,\pi )$ and belong to
$PW_{\pi -x}$ for each $x\in \lbrack d,\pi )$ as function of $\mu $.

Using the Lemma  2.1,(\ref{eq11}) gives,
\begin{equation}
\left\vert y_{R}(x,\mu )\right\vert \leq \gamma _{0}\frac{(\pi -x)}{1+|\mu
|(\pi -x)}e^{(\pi -x)|\bf{Im}\mu |}+\int_{x}^{\pi }\gamma _{0}e^{(t-x)|\bf{Im}%
\mu |}(t-x)|q(t)|\cdot |y_{R}(t,\mu )|dt  \label{eq14}
\end{equation}%
\begin{equation}
e^{-(\pi -x)|\bf{Im}\mu |}\left\vert y_{R}(x,\mu )\right\vert \leq
\gamma _{0}\frac{\pi }{1+|\mu |\pi }+\int_{x}^{\pi }\gamma _{0}\pi
|q(t)|\cdot e^{-(\pi -t)|\bf{Im}\mu |}|y_{R}(t,\mu )|dt
\label{eq15}
\end{equation}%
Gronwall's inequality (Lemma 2.2) yields,%
\begin{equation}
e^{-(\pi -x)|\bf{Im}\mu |}\left\vert y_{R}(x,\mu )\right\vert \leq
\gamma _{0}\frac{\pi }{1+|\mu |\pi }\exp \left\{ \gamma _{0}\pi
\int_{x}^{\pi }|q(t)|dt\right\}   \label{eq16}
\end{equation}%
from which we get,%
\begin{equation}
\left\vert y_{R}(x,\mu )\right\vert \leq \gamma _{0}\frac{\pi }{1+|\mu |\pi }%
\exp \left\{ \gamma _{0}\pi \int_{x}^{\pi }|q(t)|dt\right\} e^{(\pi -x)|%
\bf{Im}\mu |}\leq \frac{\gamma _{5}}{1+|\mu |\pi }e^{(\pi
-x)|\bf{Im}\mu |} \label{eq17}
\end{equation}%
where $\gamma _{5}=\gamma _{0}\pi \exp \left\{ \gamma _{0}\pi \int_{0}^{\pi
}|q(t)|dt\right\} .$

\begin{equation*}
y_{R}(x,\mu )+\frac{\sin \mu (\pi -x)}{\mu }=-\int_{x}^{\pi }\frac{\sin \mu
(x-t)}{\mu }q(t)y_{R}(t,\mu )dt
\end{equation*}%
\begin{equation}
\left\vert y_{R}(x,\mu )+\frac{\sin \mu (\pi -x)}{\mu }\right\vert \leq
\left\{ \gamma _{0}\frac{\pi }{1+|\mu |\pi }\frac{\gamma _{5}}{1+|\mu |\pi }%
\int_{0}^{\pi }|q(t)|dt\right\} e^{(\pi -x)|\bf{Im}\mu |}=\frac{\gamma _{6}%
}{\left( 1+|\mu |\pi \right) ^{2}}e^{(\pi -x)|\bf{Im}\mu |}
\label{eq18}
\end{equation}%
where $\gamma _{6}=\gamma _{0}\pi \gamma _{5}\int_{0}^{\pi }|q(t)|dt$.
Similarly,%
\begin{eqnarray}
\left\vert y_{R}^{\prime }(x,\mu )-\cos \mu (\pi -x)\right\vert
&\leq &\int_{x}^{\pi }e^{(t-x)|\bf{Im}\mu |}\left\vert
q(t)\right\vert \cdot
\left\vert y_{R}(t,\mu )\right\vert dt  \notag \\
&\leq &\left\{ \frac{\gamma _{5}}{1+|\mu |\pi }\int_{0}^{\pi }\left\vert
q(t)\right\vert dt\right\} e^{(\pi -x)|\bf{Im}\mu |}  \notag \\
&=&\frac{\gamma _{7}}{1+|\mu |\pi }e^{(\pi -x)|\bf{Im}\mu |}
\label{eq20}
\end{eqnarray}%
where $\gamma _{7}=\gamma _{5}\int_{0}^{\pi }|q(t)|dt$.

Thus, $y_{R}(x,\mu )$, $y_{R}(x,\mu )+\frac{\sin \mu (\pi -x)}{\mu }$, $%
y_{R}^{\prime }(x,\mu )-\cos \mu (\pi -x)$ are entire functions of $\mu $
for each $x\in \lbrack d,\pi )$ and belong to $PW_{\pi -x}$ for each $x\in
\lbrack d,\pi )$ as function of $\mu $.

\bigskip

Although we have obtained much higher estimates in \cite{C97} and \cite%
{C2001} to the expense of subtracting terms involving multiple
integrals, and as in \cite{C05a}, we shall stick with the
estimates given in Theorem 2.2, hence avoiding any (multiple)
integration(s) and show by the same token that we can get a higher
order estimate of the eigenvalues of the problem at hand at a very
low cost. In fact we do not have even to keep on increasing the
number of sampling points.

Let $PW_{\sigma }$~denote the Paley-Wiener space
\begin{equation*}
PW_{\sigma }=\{f\mbox{entire, }|f(\mu )|\leq C\mbox{e}^{\sigma |\mbox{Im}\mu
|}\mbox{, }\int_{\mbox{R}}|f(\mu )|^{2}d\mu <\infty \}
\end{equation*}

Let $h_{kl}$~be defined by
\begin{equation}
\left\{
\begin{array}{ccc}
h_{11}(\mu ) & = & \left( \frac{\sin \theta \mu }{\theta \mu }\right)
^{m}\left( y_{L}(d,\mu )-\cos \mu d\right)  \\
h_{12}(\mu ) & = & \left( \frac{\sin \theta \mu }{\theta \mu }\right)
^{m}\left( y_{L}^{\prime }(d,\mu )+\mu \sin \mu d\right)  \\
h_{21}(\mu ) & = & \left( \frac{\sin \theta \mu }{\theta \mu }\right)
^{m}\left( y_{R}(d,\mu )+\frac{\sin \mu (\pi -d)}{\mu }\right)  \\
h_{22}(\mu ) & = & \left( \frac{\sin \theta \mu }{\theta \mu }\right)
^{m}\left( y_{R}^{\prime }(d,\mu )-\cos \mu (\pi -d)\right)
\end{array}%
\right.   \notag
\end{equation}

Then we rewrite $y_{L}(d,\mu ),~y_{L}^{\prime }(d,\mu ),~y_{R}(d,\mu )$ and $%
y_{R}^{\prime }(d,\mu )$~as
\begin{equation}
\left\{
\begin{array}{ccc}
y_{L}(d,\mu ) & = & h_{11}(\mu )\left( \frac{\sin \theta \mu }{\theta \mu }%
\right) ^{-m}+\cos \mu d \\
y_{L}^{\prime }(d,\mu ) & = & h_{12}(\mu )\left( \frac{\sin \theta \mu }{%
\theta \mu }\right) ^{-m}-\mu \sin \mu d \\
y_{R}(d,\mu ) & = & h_{21}(\mu )\left( \frac{\sin \theta \mu }{\theta \mu }%
\right) ^{-m}-\frac{\sin \mu (\pi -d)}{\mu } \\
y_{R}^{\prime }(d,\mu ) & = & h_{22}(\mu )\left( \frac{\sin \theta \mu }{%
\theta \mu }\right) ^{-m}+\cos \mu (\pi -d)%
\end{array}%
\right.   \notag
\end{equation}

\begin{theorem}
Let $\vartheta $~be a positive constant and $m$~ be a positive integer ($%
m\geq 2$). The functions $h_{kl},~(k,l=1,2)$~ belong to the Paley space $%
PW_{\sigma }$~with $\sigma =\sigma _{0}+m\theta $ and satisfy the estimates
\begin{equation}
\left\vert h_{kl}(\mu )\right\vert \leq \frac{\gamma }{(1+\theta |\mu |)^{m}}%
e^{\sigma \left\vert \mbox{Im}\mu \right\vert }  \notag
\end{equation}%
$k,l=1,2$~for some positive constant $\gamma $, $\sigma _{0}=\max \left\{
d,\pi -d\right\} $.
\end{theorem}

\textbf{Proof:} That $h_{kl}$~are entire and satisfy the given estimates is
a direct consequence of Theorem 2.2 and the fact that $\frac{\sin \theta \mu
}{\theta \mu }$ is an entire function of $\mu $ and satisfy the estimate in
Lemma 2.1.

Since the $h_{kl}(\mu)$ belong to the Paley-Wiener space $PW_{\sigma }$ ~for
each $k,l=1,2$, they can be recovered from their values at the points\ $\mu
_{j}=j\frac{\pi }{\sigma }$, $j\in Z$, using the following celebrated
theorem,

\begin{theorem}[Whitaker-Shannon-Kotel'nikov]
Let $h\in PW_{\sigma }$, then
\begin{equation}
h(\mu )=\sum_{j=-\infty }^{\infty }h(\mu _{j})\frac{\sin \sigma (\mu -\mu
_{j})}{\sigma (\mu -\mu _{j})}  \notag
\end{equation}
$\mu _{j}=j\frac{\pi }{\sigma }$. The series converges absolutely and
uniformly on compact subsets of $C$ and in $L_{d\mu }^{2}(R)$.
\end{theorem}

For all practical purposes, we consider finite summations, therefore we need
to approximate $h_{kl}$~by a truncated series $h_{kl}^{[N]}$. The following
lemma gives an estimate for the truncation error.

\begin{lemma}[Truncation error]
Let $h_{kl}^{[N]}(\mu )=\sum_{j=-N}^{N}h_{kl}(\mu _{j})\frac{\sin \sigma
(\mu -\mu _{j})}{\sigma (\mu -\mu _{j})}$ denote the truncation of $%
h_{kl}(\mu )$. Then, for $|\mu |<N\pi /\sigma $,

\begin{displaymath}
\left| h_{kl}(\mu )-h_{kl}^{[N]}(\mu )\right|
 \leq \frac{|\sin \mu |c_{2}}{\pi (\pi /\sigma )^{m-1}\sqrt{1-4^{-m+1}}%
}\left[ \frac{1}{\sqrt{(N\pi /\sigma )-\mu }}+\frac{1}{\sqrt{(N\pi /\sigma
)+\mu }}\right] \frac{1}{(N+1)^{m-1}},
\end{displaymath}
where $c_{3}=||\mu ^{m-1}h_{kl}(\mu )||_{2}$.
\end{lemma}

\textbf{Proof:} Since $\mu ^{m-1}h_{kl}(\mu )\in L^{2}(-\infty ,\infty )$,
Jagerman's result (see \cite{Z93}, Theorem 3.21, p.90) is applicable and
yields the given estimate for the $h_{kl},~k,l=1,2$.

An approximation $B_{N}$~to the characteristic function $B$~is provided by
replacing the $h_{kl}$~by its approximation $h_{kl}^{[N]}$, and we obtain at
once,

\begin{lemma}
The approximate characteristic function $B_{N}$~satisfies the estimate,

\begin{displaymath}
\left| B(\mu )-B_{N}(\mu )\right|
 \leq \left|\frac{\sin \theta \mu }{\theta \mu }\right|^{-m}\frac{%
|\sin \mu |c_{4}}{\pi (\pi /\sigma )^{m-1}\sqrt{1-4^{-m+1}}}\left[ \frac{1}{%
\sqrt{(N\pi /\sigma )-\mu }}+\frac{1}{\sqrt{(N\pi /\sigma )+\mu }}\right]
\frac{1}{(N+1)^{m-1}},
\end{displaymath}
for some positive constant $c_{4}$.
\end{lemma}

We claim the following,

\begin{theorem}
Let $\overline{\mu}^{2}$~be an exact eigenvalue of $B$~
(multiplicity $1$) and denote by $\mu_{N}^{2}$~the corresponding
approximation of a square of a zero of $B_{N}$. Then, for
$|\mu_{N}|<N\pi/\sigma$, we have,
\begin{eqnarray*}
|\mu_{N}-\overline{\mu}|&\leq&\left(\frac{1}{|B'(\widetilde{\mu})|}%
\left|\frac{\sin \theta \mu_{N} }{\theta \mu_{N } }\right|^{-m}\frac{|\sin
\mu_{N }|c_{4}}{\pi (\pi /\sigma )^{m-1}\sqrt{1-4^{-m+1}}}\right) \\
&&\times\left[ \frac{1}{\sqrt{(N\pi /\sigma )-\mu_{N } }}+\frac{1}{\sqrt{%
(N\pi /\sigma )+\mu_{N } }}\right] \frac{1}{(N+1)^{(m-1)}}
\end{eqnarray*}
where the $\inf$~is taken over a ball centered at $\mu_{N}$~with radius $%
|\mu_{N}-\overline{\mu}|$~and not containing a multiple of $\pi/\theta$.
\end{theorem}

\textbf{Proof} Since $\overline{\mu}$~is a simple zero of $B$ , we
have,
\begin{equation}
B(\overline{\mu })-B(\mu
_{N})=(\overline{\mu}-\mu_{N})B'(\widetilde{\mu})  \notag
\end{equation}
for some $\widetilde{\mu}$. Thus,
\begin{eqnarray*}
|\overline{\mu}-\mu_{N}|&=& \frac{|B(\overline{\mu })-B(\mu
_{N})|}{
|B'(\widetilde{\mu})|} \\
&\leq&\frac{1}{|B'(\widetilde{\mu})|}\left|\frac{\sin \theta \mu_{N} }{%
\theta \mu_{N } }\right|^{-m}\frac{|\sin \mu_{N }|c_{4}}{\pi (\pi
/\sigma
)^{m-1}\sqrt{1-4^{-m+1}}} \\
&&\times\left[ \frac{1}{\sqrt{(N\pi /\sigma )-\mu_{N } }}+\frac{1}{\sqrt{%
(N\pi /\sigma )+\mu_{N } }}\right] \frac{1}{(N+1)^{m-1}}
\end{eqnarray*}
where the $\inf$~is taken over a ball centered at $\mu_{N}$~with radius $%
|\mu_{N}-\overline{\mu}|$~and not containing a multiple of $\pi/\theta$.
Thus, the result.

\section{Numerical examples}

\setcounter{equation}{0} In this section, we shall work out an
example to illustrate our method. We shall take $N=40$ , $m=6$. We
have taken $\theta
=\sigma _{0}/(N-m)$ in order to avoid the first singularity of $\left( \frac{%
\sin \theta \mu _{N}}{\theta \mu _{N}}\right) ^{-1}$. The sampling values
were obtained using the Fehlberg 4-5 order Runge-Kutta method.

\textbf{Example 3.1}
\begin{equation}
\left\{
\begin{array}{c}
-y^{\prime \prime }+q(x)y=\mu ^{2}y\qquad ,\qquad x\in (0,\pi ) \\
y^{\prime }(0)=0=y(\pi ) \\
y(d+0)=ay(d-0) \\
y^{\prime }(d+0)=a^{-1}y^{\prime }(d-0)%
\end{array}%
\right.
\end{equation}%
where $a=2$, $d=1$ and $q(x)=x$.

\begin{center}
\begin{equation*}
\begin{tabular}{|c|c|c|c|c|}
\hline
Index & Exact & RSM & Absolute Error & Relative Error \\
\hline
$1$ & $1.22788546912$ & $1.227885469249$ & $1.31357098\times 10^{-10}$ & $%
1.06978298879\times 10^{-10}$ \\
$2$ & $1.83749384727$ & $1.837493847255$ & $1.6491398678\times 10^{-11}$ & $%
8.9749408974\times 10^{-12}$ \\
$3$ & $2.68396812434$ & $2.683968124476$ & $1.32813906505\times 10^{-10}$ & $%
4.94841594058\times 10^{-11}$ \\
$4$ & $3.85661744715$ & $3.856617447367$ & $2.1294955259\times 10^{-10}$ & $%
5.52166647358\times 10^{-11}$\\
 \hline
\end{tabular}
\end{equation*}
\end{center}

The graphs of the first four exact eigenfunctions and their
approximations are displayed in Figure 1. The inner products of
any two different approximate eigenfunctions are of the order of
$10^{-10}$.
\begin{figure}[ht]
\centerline{\hbox{\psfig{figure=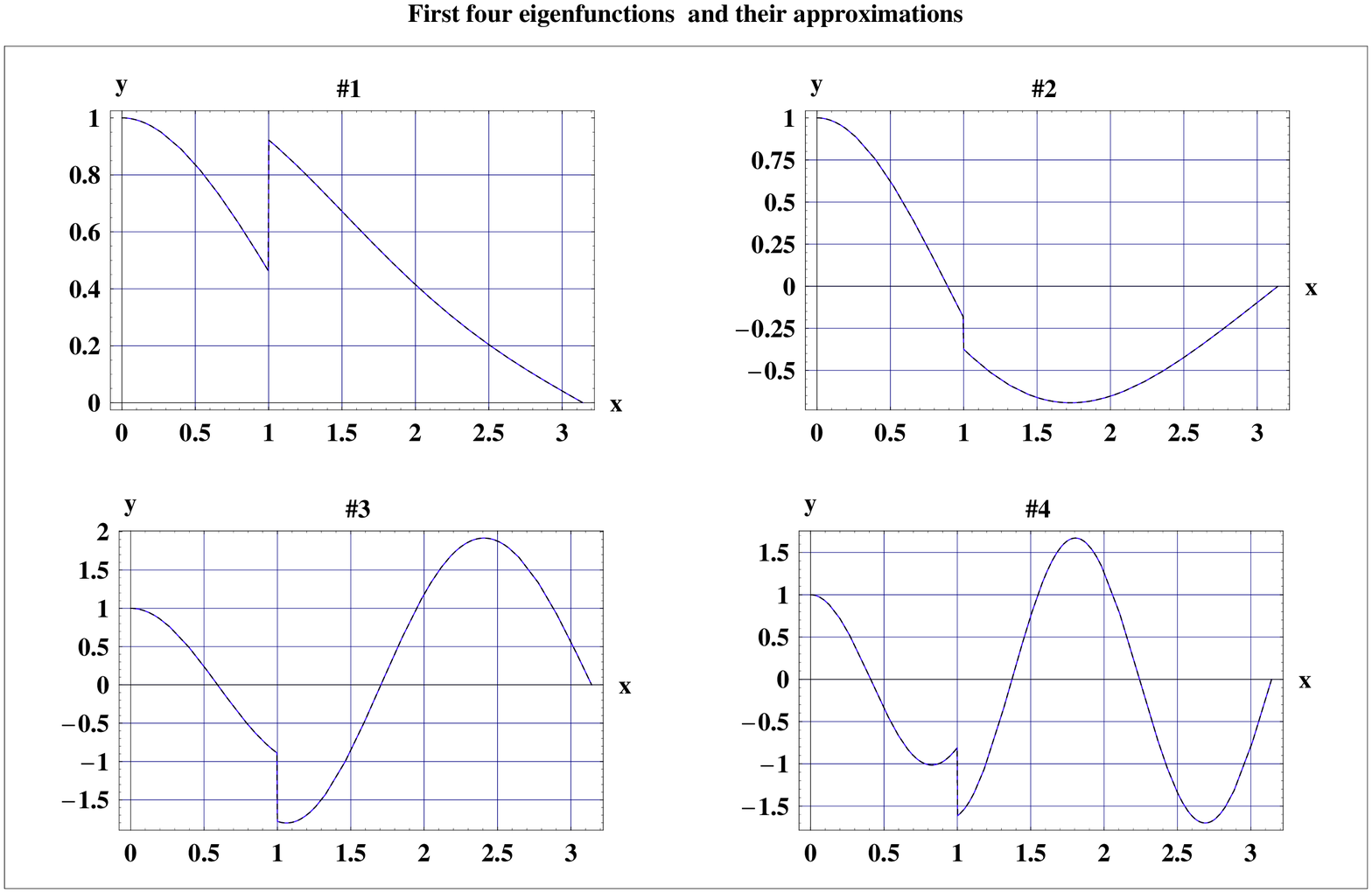,height=17cm,width=17cm}}}
\label{Figure 1}
\end{figure}

\section{Conclusion}

In this paper, we have used the \textit{regularized sampling
method} introduced recently \cite{C05a} to compute the eigenvalues
of Sturm-Liouville problems with discontinuity conditions inside a
finite interval. We recall that this method constitutes an
improvement upon the method based on Shannon sampling theory
introduced in \cite{BC96} since it uses a regularization avoiding
any multiple integration. The method allows us to get higher order
estimates of the eigenvalues at a very low cost. We have presented
an  example to illustrate the method and compared the computed
eigenvalues with the exact ones. We have also computed the
approximate eigenfunctions and compared them with the exact
eigenfunctions and checked that the inner products of any two
different approximate eigenfunctions are of the order of
$10^{-10}$, a point which validates our method. We shall present
in a future paper a generalization of the above result together
with extensive numerical computations.

\section*{Acknowledgments}

The author wishes to thank King Fahd University of Petroleum and Minerals
for its constant support.


\begin{thebibliography}{99}
\bibitem{A2006} R. Kh. Amirov, On Sturm-Liouville operators with
discontinuity conditions inside an interval, J. Math. Anal. and
Appl., 317, (2006), 163-176

\bibitem{BC96} A. Boumenir, B. Chanane, \textit{Eigenvalues of S-L systems
using sampling theory}, Applicable Analysis, Vol. 62, pp. 323-334, (1996)

\bibitem{C05a} B. Chanane, \textit{Computation of the eigenvalues
of Sturm-Liouville Problems with parameter dependent boundary
conditions using the regularized sampling method}, Math. of
Computation, 74 (2005), no. 252, 1793--1801 (published
electronically S 0025-5718(05)01717-5 in 2005).

\bibitem{C2001} B. Chanane, \textit{High Order Approximations of the
Eigenvalues of Sturm-Liouville Problems with Coupled Self-Adjoint Boundary
Conditions},Applicable Analysis, Vol. 80, pp. 317-330 (2001)

\bibitem{C97} B. Chanane, \textit{High Order Approximations of the
Eigenvalues of Regular Sturm-Liouville Problems,} J. Math. Anal. and Appl.,
226, pp.121-129, (1998)

\bibitem{DS71} N. Dunford and Schwartz, \textbf{Linear Operators, Part III},
Wiley-Interscience, New York ,(1971)

\bibitem{EE87} E. E. Edmunds and W. D. Evans, \textbf{Spectral Theory and
Differential Operators}, New York: Clarendon/Oxford University Press), (1987)

\bibitem{HS97} D. Hinton, P. W. Schaefer, \textit{Spectral Theory and
Computational Methods of Sturm-Liouville Problems}, Marcel Dekker, Inc.
(1997)

\bibitem{N68a} M. A. Naimark, \textbf{Linear Differential Operators}, Part
I, Ungar, (1968)

\bibitem{P93} J. D. Pryce, \textbf{Numerical Solution of Sturm-Liouville
Problems}, Oxford Science Publications, Clarendon Press, (1993)

\bibitem{Z93} A. I. Zayed, \textbf{Advances in Shannon's Sampling Theory},
CRC Press (1993)
\end{thebibliography}
\end{document}